\documentclass[11pt]{article}

\usepackage{latexsym,amsfonts,amsmath,amssymb,graphicx,hyperref,verbatim,enumerate}

\newcommand{\EE}{\mbox{\bf E}\,}

\newcommand{\R}{\mathbb{R}}
\newcommand{\C}{\mathbb{C}}
\newcommand{\Q}{\mathbb{Q}}
\newcommand{\HH}{\mathbb{H}}
\newcommand{\N}{\mathbb{N}}
\newcommand{\D}{\mathbb{D}}
\newcommand{\TT}{\mathbb{T}}
\newcommand{\Z}{\mathbb{Z}}

\newcommand{\pa}{\partial}

\newcommand{\F}{{\cal F}}

\newcommand{\BGE}{\begin{equation}}
\newcommand{\BGEN}{\begin{equation*}}
\newcommand{\EDE}{\end{equation}}
\newcommand{\EDEN}{\end{equation*}}

\def\til{\widetilde}

\def\sem{\setminus}

\def\lin{\overline}

 \DeclareMathOperator{\id}{id}
\DeclareMathOperator{\Imm}{Im } \DeclareMathOperator{\Ree}{Re }

 \DeclareMathOperator{\JP}{JP}

 \DeclareMathOperator{\J}{J}

\def\h0{{\bf h}}

\newtheorem{Lemma}{Lemma}[section]
\newtheorem{Theorem}{Theorem}[section]

\newtheorem{Proposition}{Proposition}[section]

\numberwithin{equation}{section}

\begin{document}
\title{\bf Loop-Erasure of Plane Brownian Motion }
\date{\empty}
\author{Dapeng Zhan}
\maketitle
\begin{abstract}
We use the coupling technique to prove that there exists a loop-erasure of a plane Brownian motion stopped on exiting a simply connected domain, and the loop-erased curve is the reversal of a radial SLE$_2$ curve.
\end{abstract}

\section{Introduction}
In \cite{Law}, the loop-erasure of a finite path on a graph is defined as follows. Let $X=(X_0,X_1,\dots,X_n)$ be such a finite path. Let $w(0)=0$ and $\tau=0$. If $X_{w(\tau)}\ne X_n$, let $w(\tau+1)=\sup\{k:X_k=X_{w(\tau)}\}+1$, $\tau=\tau+1$, and repeat this process; if $X_{w(\tau)}=X_n$, we stop. In the end, we get integer numbers $\tau\ge 0$ and $0=w(0)<w(1)<\dots<w(\tau)$. Then the lattice path $Y_k=X_{w(k)}$, $0\le k\le \tau$, is called the loop-erasure of $X$. It is easy to see that every vertex of $Y$ lies on $X$, $Y$ is a simple lattice path, and has the same initial and final vertices as $X$. 

From the definition, it is clear that a path $Y=(Y_0,Y_1,\dots,Y_\tau)$ is the loop-erasure of another path $X=(X_0,X_1,\dots,X_n)$ if and only if there is an increasing function $w:\{0,1,\dots,\tau\}\to\{0,1,\dots,n\}$ such that $w(0)=0$, $X_{w(\tau)}=X_n$, and for $0\le k\le\tau$, the path $(Y_0,\dots,Y_k)$ intersects the path $(X_{w(k)},X_{w(k)+1},\dots,X_n)$ at only one vertex, which is $Y_k=X_{w(k)}$. From this observation,
we may extend the definition of loop-erasure to (continuous) curves. Let $X(t)$, $a\le t\le b$, be a curve. We say a pair $(Y,w)$ is a loop-erasure of $X$ if $Y$ is a curve $Y(t)$, $c\le t\le d$, and $w$ is an increasing function from $[c,d]$ into $[a,b]$ such that $w(c)=a$, $X(w(d))=X(b)$, and for $c\le T\le d$, the curve $Y(t)$, $c\le t\le T$, intersects the curve $X(s)$, $w(T)\le s\le b$, at only one point, which is $Y(T)=X(w(T))$. We also say that $Y$ is a loop-erasure of $X$ if the function $w$ exists. It is easy to see that every point of $Y$ lies on $X$, $Y$ is a simple curve, and has the same initial and final points as $X$. Two loop-erasures $(Y_j,w_j)$, $j=1,2$, of $X$ are called equivalent if there is a continuous and increasing function $\theta$ that maps the domain of $Y_1$ onto that of $Y_2$ such that $w_2=w_1\circ \theta$. Given a curve $X(t)$, $a\le t\le b$, there may not exist a loop-erasure of $X$; and if a loop-erasure exists, it may not be unique up to equivalence.

In this paper we will derive the existence of a loop-erasure of a plane Brownian motion up to some finite stopping time. It is well-known that simple random walks on a regular lattice such as $\delta\Z^2$ converges to plane Brownian motions as the mesh $\delta\to0$. The loop-erasure of a simple random walk is called a loop-erased random walk (LERW). Lawler, Schramm, and Werner proved \cite{LSW-2} that the LERW on the discrete approximation of a simply connected domain converges to the Schramm-Loewner evolution (SLE) \cite{LawSLE} with parameter $\kappa=2$, i.e., SLE$_2$, when the mesh of the lattice tends to $0$. So it is reasonable to conjecture that a plane Brownian motion in a simply connected domain a.s.\ has a unique (up to equivalence) loop-erasure, which generates an SLE$_2$ curve. In this paper we will prove the existence. The uniqueness is still open to the author. In addition, we expect that there exists a deterministic algorithm to erase the loops on the Brownian motion. This is also not solved in this paper. The result in this paper extends naturally to finitely connected domains. For simplicity, we will only deal with simply connected domains.

In \cite{LSW-8/3}, Lawler, Schramm, and Werner proved that adding Brownian bubbles to an SLE$_2$ curve yields the image of a Brownian excursion in a simply connected domain. Such Brownian excursion is a Brownian motion started from a boundary point of the domain, conditioned to stay inside the domain until it leaves the domain at another boundary point. Their result gives an evidence that loop-erasures of plane Brownian motion should exist. 

We will use the coupling technique introduced in \cite{reversibility} to prove the existence of the loop-erasure. The coupling technique is used to create a coupling of a conditional plane Brownian motion with a radial SLE$_2$ curve in a simply connected domain such that, for every $t$ in the definition domain of the radial SLE$_2$ curve, say $\beta$, the first hitting point of the plane Brownian motion at the set $\beta[0,t]$ is the tip point: $\beta(t)$. Then it is easy to see that the reversal of $\beta$ is a loop-erasure of the plane Brownian motion in the domain.

\section{Plane Brownian Motion in  Simply Connected Domains}
We use the convention that a standard real Brownian motion starts from $0$, and has variance $t$ at time $t$ for $t\ge 0$, and that a standard complex Brownian motion is a complex valued random process whose real part and imaginary part are two independent standard real Brownian motions. Suppose $B_{\C}(t)$ is a standard complex Brownian motion, and $D\subsetneqq\C$ is a simply connected domain containing $0$. Let $\tau=\tau_D$ be the first time that $B_{\C}(t)\not\in D$. Then $\tau$ is an a.s.\ finite stopping time. In this paper, we will discuss loop-erasures of $B_{\C}(t)$, $0\le t\le \tau$. The following is the main theorem in this paper.

\begin{Theorem}
  Almost surely there is a loop-erasure $(\gamma,w)$ of $B_{\C}(t)$, $0\le t\le \tau$, where $\gamma$ is the reversal of a radial SLE$_2$ curve in $D$ that grows from a random prime end (\cite{Ahl}) of $D$ towards $0$, and the distribution of the prime end is the harmonic measure in $D$ seen from $0$. \label{main}
\end{Theorem}

From Riemann Mapping Theorem and conformal invariance (up to time-change) of complex Brownian motion \cite{RY}, SLE, and harmonic measure, we suffice to consider the special case that $D=\D:=\{z\in\C:|z|<1\}$. Let $P(z)=\Ree \frac{1+z}{1-z}$, $z\in\D$;
for $\rho\in\TT:=\{z\in\C:|z|=1\}$, let $P_\rho(z)=P(z/\rho)$. Then $P_\rho$ is harmonic and positive in $\D$; vanishes on $\TT$ except at $\rho$; and $P_\rho(0)=1$. We call $P_\rho$ the normalized (by its value at $0$) Poisson kernel in $\D$ with the pole at $\rho$. Let $\delta_\rho(t)$, $0\le t<\tau_\rho$, be a complex valued function that solves the ODE $\delta_\rho'(t)=\frac{2{\pa_{\lin z}} P_\rho}{P_\rho}(\delta_\rho(t)+B_{\C}(t))$ with $\delta(0)=0$, and suppose that the solution can not be extended beyond $\tau_\rho$. Here $2{\pa_{\lin z}}=\pa_x+i\pa_y$. Let $B_{\C}^\rho(t)=B_{\C}(t)+\delta_\rho(t)$, $0\le t<\tau_\rho$. Then $B_{\C}^\rho(t)$ stats from $0$, satisfies the SDE \BGE dB_{\C}^\rho(t)=dB_{\C}(t)+\frac{2{\pa_{\lin z}} P_\rho(B_{\C}^\rho)}{P_\rho(B_{\C}^\rho)}dt,\quad 0\le t<\tau_\rho;\label{SDE-BM}\EDE and there is no compact set $K\subset D$ such that $B_{\C}^\rho(t)\in K$ for $0\le t<\tau_\rho$. For otherwise, the solution $\delta_\rho(t)$ could be extended beyond $\tau_\rho$.

From It\^o's formula \cite{RY}, the process $M_\rho(t):=P_\rho(B_{\C}(t))$, $0\le t<\tau$, is a positive local martingale. So if $\sigma$ is any Jordan curve in $\D$ surrounding $0$, and $\tau_\sigma$ is the first time that $B_{\C}(t)$ visits $\sigma$, then $\EE[M_\rho(\tau_\sigma)]=M_\rho(0)=1$. From Girsanov Theorem, it is easy to check that the distribution of $(B_{\C}^\rho(t):0\le t<\tau_\sigma)$ is absolutely continuous w.r.t.\ that of $(B_{\C}(t):0\le t<\tau_\sigma)$, and the Radon-Nikodym derivative is $M_\rho(\tau_\sigma)$. Now we may decompose the distribution of $B_{\C}(t)$, $0\le t<\tau$, as follows.

\begin{Lemma}
  Let $\nu$ denote the distribution of $(B_{\C}(t):0\le t<\tau)$. For every $\rho\in\TT$, let $\mu(\rho,\cdot)$ denote the distribution of $(B_{\C}^\rho(t):0\le t<\tau_\rho)$. Then $\mu(\cdot,\cdot)$ is a random measure, and $\nu=\int_\TT \mu(\rho,\cdot)d\lambda(\rho)$, where $\lambda$ is the uniform probability measure on $\TT$. \label{random-measure}
\end{Lemma}
{\bf Proof.} We are considering probability measures on the space of curves $\gamma(t)$, $0\le t<T$, in $\D$, started from $0$. Let $(\F_t)$ denote the natural filtration generated by the curves. For each $n\in\N$, let $\tau_n$ denote the first time when $|\gamma(t)|\ge 1-1/n$. Then each $\tau_n$ is an $(\F_t)$-stopping time, and the whole sigma-algebra $\F$ is generated by the union $\bigcup_{n\in\N} \F_{\tau_n}$. From an earlier observation, for each $n\in\N$ and $\rho\in\TT$, $\mu(\rho,\cdot)$ is absolutely continuous w.r.t.\ $\nu$ on $\F_{\tau_n}$, and the Radon-Nikodym derivative is $P_\rho(B_{\C}(\tau_n))$. We have that $\rho\mapsto P_\rho(B_{\C}(\tau_n))$ is continuous, and $\int_{\TT} P_\rho(B_{\C}(\tau_n))d\lambda(\rho)=1$. Thus, $\nu=\int_\TT \mu(\rho,\cdot)d\lambda(\rho)$ on $\F_{\tau_n}$. Finally, since $\F$ is the $\sigma$-algebra generated by the union $\bigcup_{n\in\N} \F_{\tau_n}$, which is an algebra, so the proof is finished by Monotone Class Theorem. $\Box$

\vskip 3mm

Let $W_\rho(z)=\frac{\rho+z}{\rho-z}$. Then $P_\rho=\Ree W_\rho$, $W_\rho$ maps $\D$ conformally onto the right half plane $\{\Ree z>0\}$, and maps $\rho$ to $\infty$. So $P_\rho(B_{\C}^\rho(t))=\Ree Z(t)$ and $2{\pa_{\lin z}} P_\rho= \lin{W_\rho'}$. Let $Z_\rho(t)=W_\rho(B_{\C}^\rho(t))$. From It\^o's formula, $Z_\rho(t)$ satisfies the SDE: $$dZ_\rho(t)=W_\rho'(B_{\C}^\rho(t))dB_{\C}(t)+\frac{|W_\rho'(B_{\C}^\rho(t))|^2}{\Ree Z_\rho(t)}dt, \quad 0\le t<\tau_\rho. $$
Let $u_\rho(t)=\int_0^t |W_\rho'(B_{\C}^\rho(s))|^2ds$, $0\le t<\tau_\rho$. Then $u_\rho$ is continuous and increasing, and maps $[0,\tau_\rho)$ onto $[0,S_\rho)$ for some $S_\rho\in(0,\infty]$. Let $Z_\rho^u(t)=Z_\rho(u_\rho^{-1}(t))$, $0\le t<S_\rho$. Then there is another standard complex Brownian motion $\til B_{\C}(t)$ such that $Z_\rho^u(t)$ satisfies the SDE: $$dZ_\rho^u(t)= d\til B_{\C}(t)+\frac{1}{\Ree Z_\rho^u(t)}dt,\quad 0\le t<S_\rho.$$
Since the curve $B_{\C}^\rho(t)$, $0\le t<\tau_\rho$, is not contained in any compact subset of $\D$, so $Z_\rho^u(t)$, $0\le t<S_\rho$, is not contained in any compact subset of $\{\Ree z>0\}$. Thus, $S_\rho=\infty$; the real part of $Z_\rho^u$ is a completed Bessel process of dimension $3$ started from $1$; the imaginary part of $Z_\rho^u$ is a standard real Brownian motion; and the two parts are independent. So we have a.s.\ $\lim_{t\to \infty} Z_\rho^u(t)=\infty$, which implies the following lemma.

\begin{Lemma}
  Almost surely $\lim_{t\to \tau_\rho^-} B_{\C}^\rho(t)=\rho$. \label{end-pt-1}
\end{Lemma}

So we view $B_{\C}^\rho(t)$, $0\le t<\tau_\rho$, as the stopped complex Brownian motion $B_{\C}(t)$, $0\le t<\tau$, conditioned to hit $\TT$ at $\rho$.

Let $V_\rho(z)=W_\rho^{-1}(z)=\frac{z-1}{z+1}\rho$. Then $|V_\rho'(z)|=\frac{2}{|z+1|^2}$ is independent of $\rho$. For $0\le t<S_\rho=\infty$, let \BGE v_\rho(t)=\int_0^t |V_\rho'(Z_\rho^u(s))|^2ds=\int_0^t 4{|Z_\rho^u(s)+1|^{-4}}ds.\label{v-rho}\EDE Then
$$v_\rho'(u_\rho(t))=|V_\rho'(Z_\rho(t))|^2=|V_\rho'(W_\rho(B_{\C}^\rho(t)))|^2 =|W_\rho'(B_{\C}^\rho(t))|^{-2}=u_\rho'(t)^{-1}.$$
Thus, $v_\rho=u_\rho^{-1}$. From an earlier discussion, the distribution of $Z_\rho^u(t)$, $0\le t<\infty$, is independent of $\rho$. From (\ref{v-rho}) the distribution of $Z_\rho(t)=Z_\rho^u(v_\rho^{-1}(t))$, $0\le t<\tau_\rho$, is also independent of $\rho$. Since $B_{\C}^\rho(t)=V_\rho(Z_\rho(t))$, $0\le t<tau_\rho$, so we have the following lemma.

\begin{Lemma}
  For any $\rho_1,\rho_2\in\TT$, $B_{\C}^{\rho_2}(t)$, $0\le t<\tau_{\rho_2}$, has the same distribution as $R_{\rho_2/\rho_1}\circ B_{\C}^{\rho_1}(t)$, $0\le t<\tau_{\rho_1}$, where $R_{\rho_2/\rho_1}$ is a rotation: $z\mapsto z\rho_2/\rho_1$. \label{rotation}
\end{Lemma}

We know that the radial SLE in $\D$ aimed at $0$ also have rotation symmetry. From the above three lemmas, to prove Theorem\ref{main}, we suffice to show the following theorem.

\begin{Theorem}
  Almost surely there is a loop-erasure $(\gamma,w)$ of $B_{\C}^1(t)$, $0\le t\le \tau_1$, where $\gamma$ is the reversal of a radial SLE$_2$ curve in $\D$ that grows from $1$ towards $0$. \label{Thm-2}
\end{Theorem}

\section{Schramm-Loewner Evolution}
Schramm-Loewner evolution (SLE) was introduced by Oded Schramm \cite{S-SLE} to study the scaling limits of $2$-dimensional statistical lattice model at critical value, where the conformal invariance property appears in the limit. It is very successful in giving mathematical proofs of the conjectures proposed by physicists. The definition of SLE combines the Loewner's differential equation with a stochastic input. For the completeness of this paper, we now give a brief introduction of radial SLE, which is one of the major versions of SLE. The reader may refer to \cite{RS-basic} and \cite{LawSLE} for more properties of SLE.

Let $B(t)$  be a standard real Brownian motion. Let $\kappa>0$ be a parameter. Let $\xi(t)=\sqrt\kappa B(t)$, $t\ge 0$. The following differential equation is called the radial Leowner equation driven by $\xi$. \BGE \pa_t g_t(z)=g_t(z)\frac{e^{i\xi(t)}+g_t(z)}{e^{i\xi(t)}-g_t(z)},\quad g_0(z)=z.\EDE
It turns out that there is a decreasing family of domains $(D_t:0\le t<\infty)$ with $D_0=\D$ and $0\in D_t$ for all $t\ge 0$, such that each $g_t$ is defined on $D_t$, maps $D_t$ conformally onto $\D$, and satisfies $g_t(0)=0$ and $g_t'(0)=e^t$. Moreover, almost surely \BGE \beta(t):=\lim_{\D\ni z\to e^{i\xi(t)}} g_t^{-1}(z)\label{trace}\EDE exists for $0\le t<\infty$, and $\beta(t)$, $0\le t<\infty$, is a continuous curve in $\lin\D$ with $\beta(0)=1$ and $\lim_{t\to\infty}\beta(t)=0$. Such $\beta$ is called a standard radial SLE$_\kappa$ curve. If $\kappa\in(0,4]$, $\beta$ is a simple curve, intersects $\TT$ only at its initial point, and for each $t\ge 0$, $D_t=\D\sem\beta((0,t])$; if $\kappa>4$, $\beta$ is no longer a simple curve, and for each $t\ge 0$, $D_t$ is the connected component of $D_t=\D\sem\beta((0,t])$ that contains $0$. In this paper we are mostly interested in the case $\kappa=2$, so $\beta$ is a simple curve.

There is an interesting local martingale associated with radial SLE$_2$, which was used to prove the convergence of LERW to SLE$_2$ \cite{LSW-2}. Recall that $P_{e^{i\xi(t)}}$ is the normalized Poisson kernel in $\D$ with the pole at $e^{i\xi(t)}$. Since $g_t^{-1}$ maps $\D$ conformally onto $D_t=\D\sem\beta(0,t]$, fixes $0$, and has continuous extension to $\lin{\D}$, which maps $e^{i\xi(t)}$ to $\beta(t)$, so $Q_t:=P_{e^{i\xi(t)}}\circ g_t$ is the normalized (values $1$ at $0$) Poisson kernel in $D_t$ with the pole at $\beta(t)$. We have the following proposition.

\begin{Proposition}
  Let $\kappa=2$. Then for any $z\in\D$, $(Q_t(t):0\le t<T_z)$ is a local martingale, where $T_z\in(0,\infty]$ is such that $[0,T_z)$ is the maximal interval with $z\in D_t$ for $t\in[0,T_z)$ \label{prop}
\end{Proposition}

\section{Local Martingale in Two Time Variables}

Theorem \ref{Thm-2} will be proved by constructing a global commutation coupling of the process $B_{\C}^1(t)$, $0\le t<\tau_1$, with a standard radial SLE$_2$ curve $\beta(t)$, $0\le t<\infty$. The property of the global coupling will be discussed later. In this section, we will first construct a local coupling. 

First we suppose that the conditional complex Brownian motion $B_{\C}^1(t_1)$, $0\le t_1<\tau_1$, and the standard radial SLE$_2$ curve $\beta(t_2)$, $0\le t_2<\infty$, are independent. This is a trivial coupling of the above two processes.  Let $\xi(t_2)=\sqrt2 B(t_2)$ be the driving function of $\beta$, and let $g_t$ denote the radial Loewner maps. Let $(\F^1_{t_1})$ and $(\F^2_{t_2})$ be the natural filtrations generated by $B_{\C}^1(t_1)$ and $(\xi(t_2))$, respectively. Then $(\beta(t_2))$ and $(g_{t_2})$ are $(\F^2_{t_2})$-adapted. Let $${\cal D}=\{(t_1,t_2)\in[0,\tau_1)\times[0,\infty):B_{\C}^1[0,t_1]\cap \beta[0,t_2]=\emptyset\}.$$ For every $t_2\in[0,\infty)$, let ${\cal T}_1(t_2)$ be the maximal number such that $(t_1,t_2)\in\cal D$ for $t_1\in[0,{\cal T}_1(t_2))$; for every $t_1\in[0,\tau_1)$, let ${\cal T}_2(t_1)$ be the maximal number such that $(t_1,t_2)\in\cal D$ for $t_2\in[0,{\cal T}_2(t_1))$. If $\bar t_2<\infty$ is an $(\F^2_{t_2})$-stopping time, then ${\cal T}_1(\bar t_2)$ is an $(\F^1_{t_1}\times\F^2_{\bar t_2})$-stopping time; if $\bar t_1<\tau_1$ is an $(\F^1_{t_1})$-stopping time, then ${\cal T}_2(\bar t_1)$ is an $(\F^1_{\bar t_1}\times\F^2_{t_2})$-stopping time.

Let $Q_{t_2}=P_{e^{i\xi(t_2)}}\circ g_{t_2}$ be as in Proposition \ref{prop}. Since $g_0=\id$ and $\xi(0)=0$, so $Q_0(z)=P_1(z)=\frac{1+z}{1-z}$. Define $ M$ on $\cal D$ such that $$  M(t_1,t_2)=\frac{Q_{t_2}(B_{\C}^1(t_1))} {Q_0(B_{\C}^1(t_1))}.$$ It is clear that $M(t_1,0)=1$ for any $0\le t_1<\tau_1$. Since $B_{\C}^1(0)=0$ and $\Q_{t_2}(0)\equiv 1$, so $M(0,t_2)=1$ for any $0\le t_2<\infty$. 

\begin{Lemma}
  (a) For any $(\F^1_{t_1})$-stopping time $\bar t_1<\tau_1$, $M(\bar t_1,t_2)$, $0\le t_2<{\cal T}_2(\bar t_1)$, is an $(\F^1_{\bar t_1}\times\F^2_{t_2})$-local martingale. (b) For any $(\F^2_{t_2})$-stopping time $\bar t_2<\infty$, $M(t_1,\bar t_2)$, $0\le t_1<{\cal T}_1(\bar t_2)$, is an $(\F^1_{ t_1}\times\F^2_{\bar t_2})$-local martingale. \label{mart-lemma}
\end{Lemma}
{\bf Proof.} (a) This part follows immediately from Proposition \ref{prop}. \\
(b) Let $f_{\bar t_2}=Q_{\bar t_2}/Q_0$. Then $f_{\bar t_2}$ is $\F^2_{\bar t_2}$-measurable, and $M(t_1,\bar t_2)=f_{\bar t_2}(B_{\C}^1(t_1))$. Recall that $Q_0=P_1$. From (\ref{SDE-BM}) ($\rho=1$) and It\^o's formula, we see that $M(t_1,\bar t_2)$, $0\le t_1<{\cal T}_1(\bar t_2)$, satisfies the $(\F^1_{t_1}\times \F^2_{\bar t_2})$-adapted SDE:
$$d_1 M(t_1,\bar t_2)=\Ree[2\pa_z f_{\bar t_2}(B_{\C}^1(t_1))dB_{\C}(t_1)]+\Ree[2\pa_z f_{\bar t_2}(B_{\C}^1(t_1))\frac{2{\pa_{\lin z}} Q_0 (B_{\C}^1)}{Q_0(B_{\C}^1)}]dt+\frac 12\Delta f_{\bar t_2}(B_{\C}^1(t_1))dt,$$ where $2\pa_z=\pa_x-i\pa_y$, and $\Delta=\pa_x^2+\pa_y^2$. We have $f_{\bar t_2}Q_0=Q_{\bar t_2}$, and both $Q_0$ and $Q_{\bar t_2}$ are harmonic. So
$$0=\Delta Q_{\bar t_2}=4\pa_z\pa_{\lin z}(f_{\bar t_2}Q_0)=f_{\bar t_2}\Delta Q_0+Q_0\Delta f_{\bar t_2}+4\pa_z f_{\bar t_2}\pa_{\lin z}Q_0+4\pa_{\lin z}f_{\bar t_2}\pa_z Q_0$$
$$=Q_0\Delta f_{\bar t_2}+8\Ree [\pa_z f_{\bar t_2}\pa_{\lin z} Q_0].$$
So we have
\BGE d_1 M(t_1,\bar t_2)=\Ree[2\pa_z f_{\bar t_2}(B_{\C}^1(t_1))dB_{\C}(t_1)].\label{d1M}\EDE
Thus, $M(t_1,\bar t_2)$, $0\le t_1<{\cal T}_1(\bar t_2)$, is an $(\F^1_{ t_1}\times\F^2_{\bar t_2})$-local martingale. $\Box$

\vskip 3mm
Let $\J$ denote the set of Jordan curves in $\D\sem\{0\}$ that surround $0$. For every $\sigma\in\J$, let $T^1_{\sigma}$ be the first time that $B^1_{\C}(t_1)$ hits $\sigma_1$; let $T^2_{\sigma}$ be the first time that $\beta_2(t_2)$ hits $\sigma$. Then $T^j_{\sigma}$ is an $(\F^j_{t_j})$-stopping time, $j=1,2$. Let $\JP$ denote the set of $(\sigma_1,\sigma_2)\in\J^2$ such that $\sigma_1\cap\sigma_2=\emptyset$, and $\sigma_2$ surrounds $\sigma_1$. Then for any $(\sigma_1,\sigma_2)\in\JP$, $[0,T^1_{\sigma_1}]\times[0,T^2_{\sigma_2}]\subset \cal D$. 

\begin{Lemma} For any $(\sigma_1,\sigma_2)\in\JP$, $|\ln(M)|$ is bounded on $[0,T^1_{\sigma_1}]\times[0,T^2_{\sigma_2}]$ by a constant depending only on $\sigma_1$ and $\sigma_2$. \label{bdd}\end{Lemma}
{\bf Proof.} Fix $(\sigma_1,\sigma_2)\in\JP$. In this proof, a uniform constant means a constant depending only on $\sigma_1$ and $\sigma_2$; and we say a variable is uniformly bounded if its absolute value is bounded by a uniform constant. Let $N(t_1,t_2)=Q_{t_2}(B_{\C}^1(t_1))$. Since $M(t_1,t_2)=N(t_1,t_2)/N(t_1,0)$, so we suffice to show that $\ln(N)$ is uniformly bounded on $[0,T^1_{\sigma_1}]\times[0,T^2_{\sigma_2}]$. Fix $t_1\in [0,T^1_{\sigma_1}]$ and $t_2\in [0,T^2_{\sigma_2}]$. Let $D_{\sigma_j}$ denote the domain bounded by $\sigma_j$, $j=1,2$. Let $\Omega=D_{\sigma_2}\sem\lin {D_{\sigma_1}}$ and $\Omega_{t_2}=D_{t_2}\sem \lin{D_{\sigma_1}}$ for $t_2\in[0,T^2_{\sigma_2}]$. Let $m$ and $m_{t_2}$ denote the moduli of the above doubly connected domains. Then $m$ is a uniform constant. Since $\Omega$ disconnects the two boundary components of $\Omega_{t_2}$, so $m\le m_{t_2}$. Since $g_{t_2}$ maps $D_{t_2}$ conformally onto $\D$, so it maps $\Omega_{t_2}$ onto $\D\sem g_{t_2}(\lin{D_{\sigma_1}})$, which must have modulus $m_{t_2}\ge m$. Since $0\in \lin{D_{\sigma_1}}$ and $g_{t_2}(0)=0$, so $0\in g_{t_2}(\lin{D_{\sigma_1}})$. There is uniform constant $r_m\in(0,1)$ such that the modulus of $\D\sem[0,r_m]$ equals $m$. It is known that, for connected compact sets $K\subset\D$ with $0\in K$ and the modulus of $\D\sem K$ being at least $m$, the maximum of $r(K):=\sup_{z\in K}|z|$ is attained when $K=[0,r_m]$. Now $g_{t_2}(\lin{D_{\sigma_1}})$ satisfies the property of $K$, so $g_{t_2}(\lin{D_{\sigma_1}})\subset\{|z|\le r_m\}$. Since $B_{\C}^1(t_1)\in K_{\sigma_1}$, so $|g_{t_2}(B_{\C}^1(t_1))|\le r_m$. Since $N(t_1,t_2)=Q_{t_2}(B_{\C}^1(t_1))=P(g_{t_2}(B_{\C}^1(t_1))/e^{i\xi_2(t_2)})$, where $P(z)=\Ree \frac{1+z}{1-z}$, so $\frac{1-r_m}{1+r_m}\le N(t_1,t_2)\le \frac{1+r_m}{1-r_m}$. Thus, $|\ln(N)|\le \ln(\frac{1+r_m}{1-r_m})$, which is a uniform constant. $\Box$

\vskip 3mm

Now we explain the meaning of $M(t_1,t_2)$. Fix $(\sigma_1,\sigma_2)\in\JP$. Let $\mu$ denote the joint distribution of $B_\C^1(t_1)$, $0\le t_1<\tau_1$, with $\beta_2(t_2)$, $0\le t_2<\infty$, which are independent to each other. From Lemma \ref{mart-lemma} and Lemma \ref{bdd}, we have $\int M(T^1_{\sigma_1},T^2_{\sigma_2})d\mu=M(0,0)=1$. Define $\nu_{\sigma_1,\sigma_2} $ such that $d\nu_{\sigma_1,\sigma_2} /d\mu=M(T^1_{\sigma_1},T^2_{\sigma_2})$. Then $\nu_{\sigma_1,\sigma_2}$ is also a probability measure. Now suppose the joint distribution of the above two random curves is $\nu_{\sigma_1,\sigma_2}$ instead of $\mu$. Since $M=1$ when either $t_1$ or $t_2$ equals $0$, so the marginals of $\nu_{\sigma_1,\sigma_2}$ agree with those of $\mu$. Thus, $\nu_{\sigma_1,\sigma_2}$ is also a coupling measure of $B_\C^1(t_1)$, $0\le t_1<\tau_1$, with $\beta_2(t_2)$, $0\le t_2<\infty$. We now look at the behavior of the sub-curves $B_\C^1(t_1)$, $0\le t_1\le T^1_{\sigma_1}$, and $\beta_2(t_2)$, $0\le t_2\le T^2_{\sigma_2}$. Fix any $(\F^2_{t_2})$-stopping time $\bar t_2\le T^2_{\sigma_2}$. From  (\ref{SDE-BM}), (\ref{d1M}), and Girsanov Theorem, under $\nu_{\sigma_1,\sigma_2}$, there is an $(\F^1_{t_1}\times \F^2_{\bar t_2})$-standard complex Brownian motion $\til B_{\C}(t_1)$ such that $B_{\C}^1(t_1)$, $0\le t_1\le T^1_{\sigma_1}$, satisfies the $(\F^1_{ t_1}\times\F^2_{\bar t_2})$-adapted SDE:
$$dB_{\C}^1(t_1)=d\til B_{\C}(t_1)+\frac{2{\pa_{\lin z}} P_1(B_{\C}^1)}{P_1(B_{\C}^1)}dt_1+\frac{2{\pa_{\lin z}} f_{\bar t_2}(B_{\C}^1)} {f_{\bar t_2}(B_{\C}^1)}dt_1$$ \BGE =d\til B_{\C}(t_1)+\frac{2{\pa_{\lin z}} Q_{\bar t_2}(B_{\C}^1)}{Q_{\bar t_2}(B_{\C}^1)}dt_1,\label{SDE-BM-t2}\EDE
where the second equality holds because $P_1f=Q_0f=Q_{\bar t_2}$. 

\section{Coupling Measures}
Let $M$ be as in the last section. The following proposition is similar to Theorem 6.1 in \cite{reversibility}.

\begin{Proposition}
  For any $(\sigma_1^m,\sigma_2^m)\in\JP$, $1\le m\le n$, there is a continuous function $M_*(t_1,t_2)$ defined on $[0,\infty]^2$ that satisfies the following properties: \begin{itemize}
    \item[(i)] $M_*=M$ on $[0,T^1_{\sigma_1^m}]\times[0,T^2_{\sigma_2^m}]$, $1\le m\le n$;
    \item[(ii)] $M_*(t,0)=M_*(0,t)=1$ for any $t\in[0,\infty]$;
    \item[(iii)] There are constants $C_2>C_1>0$ depending only on $(\sigma_1^m,\sigma_2^m)$, $1\le m\le n$, such that $C_1\le M_*(t_1,t_2)\le C_2$ on $[0,\infty]^2$;
    \item[(iv)] For any $(\F^2_{t_2})$-stopping time $\bar t_2$, $M(t_1,\bar t_2)$, $0\le t_1\le \infty$, is an $(\F^1_{t_1}\times \F^2_{\bar t_2})$-martingale;
    \item[(v)] For any $(\F^1_{t_1})$-stopping time $\bar t_1$, $M(\bar t_1, t_2)$, $0\le t_2\le \infty$, is an $(\F^1_{\bar t_1}\times \F^2_{ t_2})$-martingale.
  \end{itemize} \label{prop-coupling}
\end{Proposition}

Let $\JP_*$ be the set of $(\sigma_1,\sigma_2)\in \JP$ such that both $\sigma_1$ and $\sigma_2$ are polygonal curves whose vertices have rational coordinates. Then $\JP_*$ is countable. Let $(\sigma_1^m,\sigma^2_m)$, $m\in\N$, be an enumeration of $\JP_*$. For each $n\in\N$, let $M_*^n$ be the $M_*$ given by the above proposition for $(\sigma_1^m,\sigma^2_m)$, $1\le m\le n$, in the above enumeration. Let the probability $\mu$ be as in the last section. For each $n\in\N$, define $\nu_n$ such that $d\nu_n=M_*^nd\mu$. From the property of $M_*$, $\int M_*^n(\infty,\infty)d\mu=M_*^n(0,0)=1$, so $\nu_n$ is a probability measure. Since $M_*^n=1$ when either $t_1$ or $t_2$ equals $0$, so $\nu_n$ is also a coupling measure of $B_\C^1(t_1)$, $0\le t_1<\tau_1$, with $\beta_2(t_2)$, $0\le t_2<\infty$. For $1\le m\le n$, since $$\EE[M_*^n(\infty,\infty)|\F^1_{T^1_{\sigma_1^m}}\times  \F^2_{T^2_{\sigma_2^m}}]=M_*^n(T^1_{\sigma_1},T^2_{\sigma_2})=M(T^1_{\sigma_1},T^2_{\sigma_2}),$$
so on $\F^1_{T^1_{\sigma_1^m}}\times  \F^2_{T^2_{\sigma_2^m}}$, $\nu_n$ equals $\nu_{\sigma_1^m,\sigma_2^m}$ defined in the last section. 

Using the argument in Section 7 in \cite{reversibility}, we can obtain a probability measure $\nu$ as a subsequential weak limit of $(\nu_n)$ in some suitable topology, which is also a coupling measure of $B_\C^1(t_1)$, $0\le t_1<\tau_1$, with $\beta_2(t_2)$, $0\le t_2<\infty$. Then for each $m\in\N$, $\nu=\nu_{\sigma_1^m,\sigma_2^m}$ on $\F^1_{T^1_{\sigma_1^m}}\times  \F^2_{T^2_{\sigma_2^m}}$. Let $\bar t_2$ be an $(\F^2_{t_2})$-stopping time with $\bar t_2 \le T^2_{\sigma_2^m}$. From the discussion at the end of the last section, we see that $B_{\C}^1(t_1)$, $0\le t_1\le T^1_{\sigma_1^m}$, satisfies (\ref{SDE-BM-t2}) for some $(\F^1_{t_1}\times \F^2_{\bar t_2})$-standard complex Brownian motion $\til B_{\C}(t_1)$.

Fix $t_2\in(0,\infty)$. For $n\in\N$, define $$R_n=\sup\{T^1_{\sigma_1^m}:1\le m\le n, t_2\le T^2_{\sigma_2^m}\}.$$ Fix $n\in\N$. Then for any $1\le m\le n$, if $t_2\le T^2_{\sigma^2_m}$, then $B_{\C}^1(t_1)$, $0\le t_1\le T^1_{\sigma_1^m}$, satisfies (\ref{SDE-BM-t2}). So $B_{\C}^1(t_1)$, $0\le t_1\le R_n$, should also satisfy (\ref{SDE-BM-t2}).

From the definition, ${\cal T}_1(t_2)$ is the maximal number such that $B_{\C}^1(t_1)$ is disjoint from $\beta[0,t_2]$ for $0\le t_1<{\cal T}_1(t_2)$. It is easy to check that ${\cal T}_1(t_2)=\sup_{n\in\N} R_n$. Thus, $B_{\C}^1(t_1)$, $0\le t_1<{\cal T}_1(t_2)$, should also satisfy (\ref{SDE-BM-t2}). Let $W_{t_2}(z)= \frac{e^{i\xi(t_2)}+g_{t_2}(z)}{e^{i\xi(t_2)}-g_{t_2}(z)}$. Then $Q_{t_2}=\Ree W_{t_2}$; $W_{t_2}$ maps $D_{t_2}$ conformally onto the right half plane, and maps $\beta(t_2)$ to $\infty$. The argument in the proof of Lemma \ref{end-pt-1} can be used here to show that a.s.\ $\lim_{t_1\to {\cal T}_1(t_2)} B_{\C}^1(t_1)=\beta(t_1)$. Thus, $B_{\C}^1({\cal T}_1(t_2))=\beta(t_1)$. In fact, we may view $B_{\C}^1(t_1)$, $0\le t_1<{\cal T}_1(t_2)$, as the complex Brownian motion $B_{\C}(t)$ conditioned to leave $D_{t_2}$ at $\beta(t_2)$. Since this result holds for every $t_2\in(0,\infty)$. So a.s.\ for every $t_2\in\Q\cap(0,\infty)$, we have $B_{\C}^1({\cal T}_1(t_2))=\beta(t_2)$. 

From the definition, it is clear that ${\cal T}_1$ is a decreasing function, and $B_{\C}^1[0,{\cal T}_1(t_2))$ is disjoint from $\beta[0,t_2]$ for any $t_2\in[0,\infty)$. For any $a\in\R$, it is easy to check that $\{t_2:{\cal T}_1(t_2)>a\}$ is an open subset of $[0,\infty)$. So ${\cal T}_1$ is right-continuous. Since both $B_{\C}^1$ and $\beta$ are continuous, and $\Q\cap(0,\infty)$ is dense in $(0,\infty)$, so a.s.\ for any $t_2\in(0,\infty)$, $B_{\C}^1({\cal T}_1(t_2))=\beta(t_2)$. We may define $B_{\C}^1(\tau_1)=1$, $\beta(\infty)=0$, and ${\cal T}_1(\infty)=0$. Then a.s.\  $B_{\C}^1({\cal T}_1(t_2))=\beta(t_2)$ holds for any $t_2\in[0,\infty]$. Let $\gamma(t)=\beta(1/t)$, $0\le t\le \infty$, and $w(t)={\cal T}_1(1/t)$. Then $\gamma$ is a reversal of $\beta$, and a.s.\ $(\gamma,w)$ is a loop-erasure of $B_{\C}^1(t_1)$, $0\le t_1\le \tau_1$. So we finish the proof of Theorem \ref{Thm-2}. 

\section{Some Remarks}
\begin{itemize}
  \item One can prove that, under the new coupling measure $\nu$, for any $(\F^2_{t_2})$-stopping time $\bar t_2<\tau_1$, the curve $\beta(t_2)$, $0\le t_2<{\cal T}_2(\bar t_1)$, is a stopped radial SLE$_2$ curve in $\D$ started from $1$ aimed at $B_{\C}^1(\bar t_2)$. 
  \item Up to a time-change, the curve $\gamma$, i.e., the reversal of $\beta$, has the distribution of a disc SLE$_2$ \cite{Zhan} curve in $\D$ that grows from $0$ to $1$. This is a special case of the continuous LERW started from interior points defined in \cite{int-LERW}.
  \item Theorem \ref{main} can be extended to finitely connected plane domains. Let $B_{\C}(t)$, $0\le t\le \tau$, be a complex Brownian motion started from an interior point $z_0$ in a finitely connected domain $D$, stopped on hitting $\pa D$. Then it has a loop-erasure $(\gamma,w)$, where $\gamma$ is a continuous LERW from $z_0$ to a random prime end of $D$, whose distribution is the harmonic measure in $D$ viewed from $z_0$. We may also derive the loop-erasure of the complex Brownian motion in $D$ conditioned to hit a marked interior point $z_1$ in $D$, which satisfies the SDE (\ref{SDE-BM}) for $B_{\C}^\rho(t)$ with $P_\rho$ replaced by $G_D(z_1,\cdot)$, where $G_D(\cdot,\cdot)$ is the Green function in $D$. The loop-erased curve is a continuous LERW in $D$ from $z_0$ to $z_1$.
  \item We may also derive the existence of a loop-erasure of a Brownian excursion in a finitely connected domain. The loop-erased curve is a continuous LERW started from boundary \cite{LERW}. For the proof, we may use the coupling technique to construct a coupling in the degenerate case, see Section 4.4 in \cite{duality}. For example, if $X(t)$ is a standard real Brownian motion, $Y(t)$ is a Bessel process of dimension $3$ started from $0$, and $(X)$ is independent of $(Y)$, then $Z(t)=X(t)+iY(t)$ is a Brownian excursion in the upper half plan $\HH:=\{z\in\C:\Imm z>0\}$ that grows from $0$ to $\infty$. We can conclude that it has a loop-erasure, and the loop-erased curve is the continuous LERW in $\HH$ from $0$ to $\infty$, which is just the standard chordal SLE$_2$ curve.
\end{itemize}

\end{document}